\documentclass{article}

\usepackage[skip=2pt,font=scriptsize]{caption}
\usepackage{mathtools,comment}
\usepackage{accents}
\usepackage{caption,xcolor}
\usepackage{subcaption}
\usepackage{pgfplots}
\pgfplotsset{compat=1.13}
\usepgfplotslibrary{fillbetween}
\usepackage{tikz}
\newlength\figureheight
\newlength\figurewidth
\usepackage{tikz-qtree}
\usetikzlibrary{trees}
\usetikzlibrary{shapes,arrows}
\usetikzlibrary{fit,calc,arrows,positioning}                                                        
\usepackage{balance}
\usepackage{amssymb}
\usepackage{amsmath} 
\usepackage{algorithm}
\usepackage{algorithmic}
\usepackage{graphicx,psfrag}
\usepackage{booktabs}
\usepackage{amsfonts}
\usepackage{amsmath,amsthm}
\usepackage{balance}
\usepackage{algorithm}
\usepackage{algorithmic}
\usepackage{tikz}
\usepackage{tikz-qtree}
\usetikzlibrary{trees}
\usetikzlibrary{fit,calc,arrows,positioning}
\usepackage{pgfplots}
\pgfplotsset{compat=newest}
\pgfplotsset{every axis/.append style={every x tick label/.append style={font=\tiny,yshift=-0.7ex,rotate=0},
                    every y tick label/.append style={font=\tiny,xshift=-0.7ex},ylabel shift = -10 pt,
xlabel shift = -6 pt,xlabel={\tiny{$x_{p1}$}},          
                    ylabel={\tiny{$x_{p2}$}}, title style={yshift=-1.5ex}, title style={font=\tiny}}}
\pgfplotsset{plot coordinates/math parser=false}
\usepackage{caption}
\usepackage{mathtools}
\usetikzlibrary{shapes,arrows}
\def\BState{\State\hskip-\ALG@thistlm}

\newcommand{\R}{\mathbb{R}}

\newtheorem{theorem}{Theorem}

\newtheorem{corollary}{Corollary}

\newtheorem{definition}{Definition}

\newtheorem{proposition}[theorem]{Proposition}
\newtheorem{lemma}{Lemma}

\newtheorem{fact}{Fact}

\begin{document}

\title{\LARGE \bf  Parameter identification for an uncertain reaction-diffusion equation via setpoint regulation}

\author{Gildas Besançon\thanks{G. Besançon is with Univ. Grenoble Alpes, CNRS, Grenoble INP, GIPSA-lab, Grenoble, France, email: \texttt{gildas.besancon@grenoble-inp.fr}}, Andrea Cristofaro\thanks{A. Cristofaro is with Department of Computer, Control and Management Engineering, Sapienza University of Rome, Italy, email: \texttt{andrea.cristofaro@uniroma1.it}}, Francesco Ferrante\thanks{F. Ferrante is with Department of Engineering, University of Perugia, Italy, email: \texttt{francesco.ferrante@unipg.it}}
}%

\maketitle
\begin{abstract}
The problem of estimating the reaction coefficient of a system governed by a reaction-diffusion partial differential equation is tackled. An estimator relying on boundary measurements only is proposed. The estimator is based upon a setpoint regulation strategy and leads to an asymptotically converging estimate of the unknown reaction coefficient. The proposed estimator is combined with a state observer and shown to provide an asymptotic estimate of the actual system state. A numerical example supports and illustrates the theoretical results.
\end{abstract}
\section{Introduction}
Reaction-diffusion processes can be found in various areas, mostly arising from chemistry, but also from physics, biology or ecology. They are described by parabolic partial differential equations, including both diffusive and reactive phenomena. In the case of linear reaction, the model is thus characterized by two parameters: diffusion coefficient and reaction rate. In this paper, we are interested in the case of unknown reaction coefficient. Within this setting, we consider a scenario in which measurements can be collected only at one boundary. At the same boundary, we assume that the spatial derivative of the system state is controlled (Neuman's condition). The problem of adaptive control design for such a class of systems has been deeply studied (see e.g. \cite{Krstic2008}, \cite{guo2020robust},  \cite{guo2020adaptive}, \cite{smyshlyaev2007adaptive}, \cite{jian2012adaptive}, \cite{li2020adaptive}), but with full state knowledge or via the use of backstepping transformations and in a control design perspective. Recently, in \cite{ahmed2023finite} the problem of estimating the state of a reaction-diffusion equation in the presence of unknown parameters in the output equation has been addressed via a finite-dimensional adaptive observer. 

In the present paper, we pursue a different approach. In particular, we assume that the system state is not known and we use the system control input (setpoint regulation) for an identification purpose (reaction coefficient estimation). This in turns allows to provide an estimator for the full state. More precisely, the key idea consists of designing a dynamic, finite-dimensional, compensator that drives the boundary output towards a reference value. Accordingly, the unknown reaction coefficient is recovered by relying on the nonlinear mapping between the compensator steady-state and the corresponding reference value, which turns out to be a static function of the reaction coefficient. In particular, owing to the convergence of the compensator state to its steady-state value, the sought reaction coefficient can be asymptotically estimated by relying on the compensator state. 
The proposed scheme is shown to be exploitable in combination with a classical infinite-dimensional observer, thereby providing the simultaneous estimation of the unknown parameter and of the actual system state. Summarizing, the paper's contribution is twofold:
\begin{itemize}
\item We deliver an adaptive control strategy for estimating the unknown reaction coefficient;
\item We show how to apply the previous result to the design of an adaptive PDE observer.
\end{itemize}
Related results can be found in \cite{feng2017new}, 
 where unknown input observers are considered for parabolic equations with additive uncertain terms, \cite{ascencio2016adaptive} and \cite{ahmed2015adaptive}, where backstepping-based observer design is considered, or \cite{cristofaro2023adaptive}, where adaptive multiple-models are proposed to approximate the behavior of the uncertain system. 
 
The remainder of the paper is organized as follows. The problem under consideration is more formally stated in Section \ref{sec:pb}, and the main parameter estimation result presented in Section \ref{sec:main}. An extension to joint state and parameter estimation is then discussed in Section \ref{sec:joint}, while an illustrative simulation result is provided in Section \ref{sec:ex}. 
\subsection{Notation}
The notation $\R$ indicates the set of real numbers. The symbol $D(A)$ stands for the domain of the operator $A$. Given a function $f\colon U\to V$, with $U$ and $V$ being normed linear vector spaces, $df(x)$ stands (when it exists) for the Fréchet differential of $f$ at $x\in U$. The symbol $\mathcal{L}^2(0,1)$ stands for the collection of equivalence classes of measurable functions $f\colon [0,1]\to \R$ that are square integrable and $\mathcal{H}^2(0,1)\coloneqq \{f\in\mathcal{L}_2(0,1)\colon f, f^\prime\,\,\text{are\,\,absolutely\,\,continuous\,\,and}\,f^{\prime\prime}\in\mathcal{L}_2(0,1)\}$, where\footnote{An element of $\mathcal{L}_2(0,1)$ is said to be absolutely continuous if there is an absolutely continuous function in the equivalence class of the element itself.} $f^\prime$ and $f^{\prime\prime}$ stand, respectively, for the first and second (weak) derivative of the function $f$.
\subsection{Preliminaries}
In this paper, we consider semilinear abstract dynamical systems of the form
\begin{equation}
\label{eq:gen_abstract}
\dot{x}=Ax+f(x)    
\end{equation}
where $x\in\mathcal{X}$ is the system state, $\mathcal{X}$ is a real Hilbert space equipped with the natural norm $\Vert \cdot\Vert$ induced by the inner product, $A\colon D(A)\subset\mathcal{X}\to\mathcal{X}$ is a linear (unbounded) operator, and $f\colon\mathcal{X}\to\mathcal{X}$ is such that $f(0)=0$. A solution to \eqref{eq:gen_abstract} is any continuous function $\varphi\colon [0, T)\to\mathcal{X}$, with $T\in \mathbb{R}_{\geq 0}\cup\{\infty\}$, such that for all $t\in[0, T)$
$$
\begin{aligned}
&\int_0^t\varphi(s)ds\in D(A)\\
&\varphi(t)=\varphi(0)+A\int_0^t\varphi(s)ds+\int_0^t f(\varphi(s))ds. 
\end{aligned}
$$
In particular, we say that $\varphi$ is maximal if its domain cannot be extended and complete if its domain is unbounded. 

The following notion of well-posedness is considered throughout the paper.
\begin{definition}
\label{def:well-posed}
System \eqref{eq:gen_abstract} is said to be well posed if the following items are satisfied:
\begin{itemize}
    \item[$(i)$] The operator $A$ generates a strongly continuous semigroup on the Hilbert space $\mathcal{X}$;
    \item[$(ii)$] The function $f$ is locally Lipschitz continuous on the Hilbert space $\mathcal{X}$.
\end{itemize}
\end{definition}
The satisfaction of the two items outlined in Definition~\ref{def:well-posed} ensures that for all $\varphi_0\in\mathcal{X}$, there exists a unique maximal solution $\varphi$ to \eqref{eq:gen_abstract} such that $\varphi(0)=\varphi_0$. Moreover, if $\varphi$ is bounded, then it is complete; see \cite[Chapter 6.1]{pazy2012semigroups} for more details about these aspects. 

The following notions of stability are used in the paper. 
\begin{definition}
\label{def:stab}
Let items $(i)-(ii)$ in Definition~\ref{def:well-posed} hold. The origin of \eqref{eq:gen_abstract} is said to be\footnote{In light of the inherent boundedness requirement of the stability, global uniform attractivity, and local exponential stability notions, under the satisfaction of items $(i)-(ii)$ in Definition~\ref{def:well-posed}, maximal solutions to \eqref{eq:gen_abstract} turn out to be complete. This is why completeness of maximal solutions is tacitly assumed throughout Definition~\ref{def:stab}.} :
\begin{itemize}
    \item Stable if for all $\varepsilon>0$, there exists $\delta>0$ such that any maximal solution $\varphi$ to  \eqref{eq:gen_abstract} with $\Vert\varphi(0)\Vert\leq \delta$ satisfies $\Vert \varphi(t)\Vert\leq\varepsilon$ for all $t\in\R_{\geq 0}$;
    \item Globally uniformly attractive if every solution to \eqref{eq:gen_abstract} is bounded and for any $\varepsilon,\mu>0$ there exists $T>0$ such that $t\geq T$ implies $\Vert\varphi(t)\Vert\leq\varepsilon$ for every solution $\varphi$ to \eqref{eq:gen_abstract} with $\Vert\varphi(0)\Vert\leq\mu$;
    \item Globally asymptotically stable (\emph{GAS}) if both stable and globally uniformly attractive;
    \item Locally exponentially stable (\emph{LES}) if there exist $r, \lambda, \kappa>0$ such that any maximal solution $\varphi$ to \eqref{eq:gen_abstract} with $\Vert\varphi(0)\Vert\leq r$ satisfies
    $$
    \Vert \varphi(t)\Vert\leq \kappa e^{-\lambda t}\Vert\varphi(0)\Vert\quad\forall t\geq 0.
    $$
\end{itemize}
\end{definition}
\section{Problem setting}\label{sec:pb}
Let us consider the reaction-diffusion problem
\begin{equation}\label{eq:sys}
\begin{array}{rcl}
u_t(x,t)&=&\lambda u_{xx}(x,t)-k u(x,t),\,\, x\in[0,1],\, t\geq 0\\
u_x(0,t)&=&0\\
u_x(1,t)&=&v(t)\\
u(x,0)&=&u_0(x)\\
y(t)&=& u(1,t)
\end{array}
\end{equation}
for a control input function $v$ and some unknown initial state function $u_0$, with known coefficient $\lambda>0$, and unknown $k>0$, and available measurement $y$. The goal of this paper can be outlined as:
\begin{itemize}
\item Recover the value of the reaction coefficient $k$;
\item Design an asymptotic state observer for the equation~\eqref{eq:sys}.
\end{itemize}

The first intuitive attempt to tackle the estimation of \eqref{eq:sys} might be using a classical adaptive control strategy.
In particular, one may consider an adaptive observer $(\hat{u},\hat{k})$, with infinite-dimensional part to be updated by replicating the system dynamics and inserting an output injection with gain $
\varrho>0$. 
A good adaptation law for $\hat{k}$ should be found instead by looking at the derivative of a Lyapunov function $V(\cdot)$ for the error system $e(x,t)=u(x,t)-\hat{u}(x,t)$, searching among the solutions to the inequality $\dot{V}(\cdot)\leq 0$, and finally conclude by LaSalle's invariance principle. Developing the latter idea, based on the ensuing Lyapunov analysis, it becomes clear that one would need an adaptation law capable of compensating a term in a form similar to\footnote{This particular form arises when using $V(e)=\int_0^1e(x)^2dx$, which is the simplest and most natural choice for the Lyapunov function.}
$$
(k-\hat{k})\int_0^1e(x,t)\hat{u}(x,t)dx.
$$
However, this is not feasible relying on the output $y(t)$ only. In fact, the integral term above contains the whole error state $e(x,t)$ and therefore it is not available for measurement, thus hampering the chance of successfully estimating the state $u(x,t)$ and the constant $k$ by following this approach. 

The unsuccessful attempt using the classical approach motivates the search for alternative methods. We propose here a novel approach, based on a particular choice for the control input, which aims at setpoint regulation of the boundary output.

\section{Main results}\label{sec:main}
\subsection{Outline of the proposed estimator}
A building block of the proposed estimator is given by the following (finite-dimensional) dynamic feedback compensator  
\begin{equation}
\label{eq:compensator}
\begin{array}{rcl}
\dot{\chi}&=&-\rho(y-y_r) \\
v&=&\chi-\gamma (y-y_r)
\end{array}
\end{equation}
where $\gamma, \rho>0$ and $y_r\in\R$ are tuning parameters that are selected later. The closed-loop system becomes
\begin{equation}
\label{eq:PI_interc}
\begin{array}{rcl}
u_t(x,t)&=&\lambda u_{xx}(x,t)-k u(x,t),\,\, x\in[0,1],\, t\geq 0\\
\dot{\chi}&=&-\rho(u(1,t)-y_r) \\
u_x(0,t)&=&0;\,u(x,0)\,=\, u_0(x)\\
u_x(1,t)&=&\chi-\gamma (u(1,t)-y_r)\\
y(t)&=& u(1,t)
\end{array}
\end{equation}
For any $y_r\in\R$, system \eqref{eq:PI_interc} admits a unique equilibrium $\bar{z}\coloneqq (\bar{u}, \bar{\chi})\in\mathcal{H}^2(0,1)\oplus\R$. In particular, such an equilibrium corresponds to the unique solution to the following boundary value problem
\begin{equation}
\label{eq:PI_eq}
\begin{array}{rcl}
0&=&\lambda \bar{u}_{xx}(x)-k\bar{u}(x),\,\, x\in(0,1)\\
0&=&\bar{u}(1)-y_r\\
\bar{u}_x(0)&=&0\\
\bar{u}_x(1)&=&\bar{\chi}.
\end{array}
\end{equation}
By solving \eqref{eq:PI_eq} one gets
$$
\begin{aligned}
&\bar{u}(x)=\frac{\cosh\left(\sqrt{\frac{k}{\lambda}}x\right)}{\cosh\left(\sqrt{\frac{k}{\lambda}}\right)}y_r,&x\in[0,1]\\
& \bar{\chi}= \sqrt{\frac{k}{\lambda}} \tanh\left(\sqrt{\frac{k}{\lambda}}\right)y_r.
\end{aligned}
$$
To make the proposed parameter estimation strategy effective, one needs to establish conditions on the parameter $\gamma$ ensuring convergence of solutions to \eqref{eq:PI_interc} to the equilibrium $\bar{z}$. To this end, it is convenient to define the following error coordinates $\tilde{u}\coloneqq u-\bar{u}$ and $\tilde{\chi}\coloneqq \chi-\bar{\chi}$ and analyze the corresponding dynamics. In particular, by using \eqref{eq:PI_eq}, one gets:  
\begin{equation}
\begin{array}{rcl}
\tilde u_t&=&\lambda \tilde u_{xx}-k\tilde u\\
\tilde u_x(0)&=&0\\
\tilde u_x(1)&=&-\gamma \tilde u(1) +\tilde{\chi}\\
\dot{\tilde{\chi}}&=&-\rho  \tilde u(1).
\end{array}
\label{eq:error_dyns}
\end{equation}
By taking as a state $\tilde{z}\coloneqq (\tilde{u},\tilde{\chi})\in\mathcal{X}\coloneqq \mathcal{L}^2(0,1)\oplus\R$, the dynamics in \eqref{eq:error_dyns} can be formally modeled by the following abstract dynamical system
\begin{equation}
\label{eq:error_dyn}
\dot{\tilde{z}}=A_e\tilde{z}
\end{equation}
where
\begin{equation}
\begin{aligned}
&A_e\begin{bmatrix}
\tilde{u}\\
\tilde{\chi}
\end{bmatrix}\coloneqq\begin{bmatrix}
\lambda \tilde u_{xx}-k\tilde u\\
-\rho \tilde u(1)
\end{bmatrix}\\
&D(A_e)\coloneqq \left\{(\tilde{u}, \tilde{\chi})\in\mathcal{H}^2(0,1)\oplus\R\colon\tilde u_x(0)=0, \right.\\
&\hspace{5.1cm}\left.\tilde u_x(1)=-\gamma \tilde u(1) +\tilde{\chi} \right\}.\\
\end{aligned}
\label{eq:Ae_op}
\end{equation}
\subsection{Well-posedness}\label{sec:wellp}
Next, well-posedness of \eqref{eq:error_dyn} is established. 
\begin{proposition}
\label{prop:wellposed}
The unbounded operator $A_e$ defined in \eqref{eq:Ae_op} generates a strongly continuous semigroup on the Hilbert space $\mathcal{X}$ equipped with the following inner product
$$
\langle (u_1,\chi_1), (u_2,\chi_2)\rangle_{\mathcal{X}}\coloneqq \langle u_1, u_2\rangle_{\mathcal{L}^2(0,1)}+\frac{\lambda}{\rho} \chi_1 \chi_2.
$$
\end{proposition}
\begin{proof}
Observe that 
$$
A_e=A_0+T
$$
with
$$
\begin{aligned}
&T\begin{bmatrix}
\tilde{u}\\
\tilde{\chi}
\end{bmatrix}\coloneqq\begin{bmatrix}
-k\tilde u\\
0
\end{bmatrix}, \quad A_0\begin{bmatrix}
\tilde{u}\\
\tilde{\chi}
\end{bmatrix}\coloneqq\begin{bmatrix}
\lambda \tilde u_{xx}\\
-\rho \tilde u(1)
\end{bmatrix}\\
&D(T)=\mathcal{X},\quad \quad D(A_0)=D(A_e)
 \end{aligned}
$$
Thus, since $T$ is bounded, from \cite[Theorem 3.2.1, page 110]{curtain2012introduction} to show the result it is enough to show that $A_0$ generates a strongly continuous semigroup on the Hilbert space $\mathcal{X}$. To establish this property, we make use of the Lumer-Phillips's theorem \cite[Theorem 4.3, page 14]{pazy2012semigroups}. As a first step, we show that the operator 
$
A_0-\lambda I
$
is surjective. To this end, let $f=(f_1, f_2)\in\mathcal{X}$. We show that there exists $z=(u,\chi)\in D(A_0)$ such that
$$
(A_0-\lambda I)z=f.
$$
From the definition of $A_0$, the above identity gives
\begin{equation}
\begin{aligned}
&u_{xx}-u=\frac{1}{\lambda}f_1\\
&u_x(0)=0\\
&u_x(1)+\gamma u(1)-\chi=0\\
&-\rho u(1)-\lambda\chi=f_2.
\end{aligned}
\label{eq:surj_sys}
\end{equation}
Therefore, to show the above surjectivity property we show that the above system of equations admits a solution. 
By standard integration theory and using the fact that $u_x(0)=0$, it turns out that
$$
u(x)=\kappa\cosh(x)+\int_0^x \cosh(x-s)f_1(s)ds
$$
where $\kappa$ is an integration constant to be determined. At this stage, the solvability of \eqref{eq:surj_sys} with respect to $u$ and $\chi$ reduces to the solvability of the following system of equations in the variables $\kappa$ and $\chi$
\begin{equation}
\begin{aligned}
&u_x(1)+\gamma u(1)-\chi=0\\
&-\rho u(1)-\lambda\chi=f_2.
\end{aligned}
\label{eq:surj_sys2}
\end{equation}
Observe that
$
u_x(x)=\kappa\sinh(x)+\int_0^x \sinh(x-s)f_1(s)ds
$
Hence, \eqref{eq:surj_sys2} reduces to
\begin{equation}
\begin{aligned}
&\chi=\frac{-f_2-\gamma u(1)}{\lambda}\\
&\kappa\varpi=f_2-\int_0^1 (\sinh(1-s)+\cosh(1-s))f_1(s)ds
\end{aligned}
\label{eq:surj_sys3}
\end{equation}
where $\varpi\coloneqq\left(\gamma\frac{1+\lambda}{\lambda}\cosh(1)+\sinh(1)\right)$. Hence, since $\varpi$ is nonzero, the second equality in \eqref{eq:surj_sys3} is solvable. This implies that \eqref{eq:surj_sys} admits a unique solution for any $f$, thereby showing the sought surjectivity property. We now conclude the proof of the result by showing that 
\begin{equation}
\label{eq:dissy}
\langle Az, z\rangle_\mathcal{X}\leq 0\quad\forall z\in D(A_0).
\end{equation}
The latter in turn implies that $A_0$ is a dissipative operator, i.e, $\langle (\varpi I-A_0)z,  (\varpi I-A_0)z\rangle_{\mathcal{X}}\geq\varpi \langle z, z\rangle_\mathcal{X}$ for all $\varpi >0, z\in D(A_0)$. For every $z\in D(A_0)$ one has
$$
\langle Az, z\rangle_\mathcal{X}=\lambda \int_0^1 u  u_{xx}dx-\lambda\chi u(1).
$$
Therefore, by integrating by parts
$$
\begin{aligned}
\langle Az, z\rangle_\mathcal{X}=&\lambda\left(u(1) u_x(1)-u(0)u_x(0)-\int_{0}^1 u_x^2 dx\right)\\
&-\lambda\chi u(1).
\end{aligned}
$$
By recalling that $z\in D(A_0)$, one has that $u_x(0)=0, u_x(1)=-\gamma u(1)+\chi$, hence
$$
\langle Az, z\rangle_\mathcal{X}=-\lambda\int_{0}^1 u_x^2 dx-\gamma u(1)^2
$$
which implies \eqref{eq:dissy}. This concludes the proof.
\end{proof}
\subsection{Convergence analysis}\label{sec:conver}
Let us now address the convergence of the error system. The following property holds for system \eqref{eq:error_dyn}.
\begin{proposition}\label{prop:GES}
Let $\gamma>0$. Then, the origin of \eqref{eq:error_dyn} is GES.
\end{proposition}
\begin{proof}
Let for all $z\in\mathcal{X}$
\begin{equation}\label{eq:Lyapunov}
V(\tilde{z})\coloneqq \int_0^1\begin{bmatrix}
  \tilde{u}(x)\\
  \tilde{\chi}
\end{bmatrix}^\top Q\begin{bmatrix}
  \tilde{u}(x)\\
  \tilde{\chi}
\end{bmatrix}dx,
\end{equation}
where
$$
Q\coloneqq\begin{bmatrix}
   1&-\varepsilon \\
    -\varepsilon &\frac{\lambda}{\rho}
\end{bmatrix}
$$
with $\epsilon\in(0,\sqrt{\lambda/\rho})$ to be selected later.
Observe that the above selection of $\varepsilon$ ensures that $Q\succ 0$, which in turn implies that
\begin{equation}
\label{eq:sandwich}
c_1\langle z, z\rangle_{\mathcal{X}}^2 \leq V(z)\leq c_2\langle z, z\rangle_{\mathcal{X}}^2\quad \forall z\in\mathcal{X}
\end{equation}
with 
$$
c_2\coloneqq \lambda_{\max}(Q)\frac{\rho}{\lambda}, c_1\coloneqq \lambda_{\min}(Q)\frac{\rho}{\lambda}. 
$$
Let for all $\tilde{z}\in D(\mathcal{A}_{cl})$  
$$
\dot{V}(\tilde{z})\coloneqq dV(\tilde{z})\mathcal{A}_{e}\tilde{z}.
$$
Then, with some computations one gets
\begin{equation}
    \label{eq:Vdot}
\begin{aligned}
 \dot{V}(\tilde{z})=&-2k \int_0^1\tilde{u}^2+ 2\lambda\underbrace{\int_0^1  \tilde{u}\tilde{u}_{xx}}_{\omega}-\varepsilon \int_0^1 (\lambda \tilde{u}_{xx}-k\tilde{u})\tilde{\chi}\\
 &+\varepsilon\rho \int_0^1 \tilde{u}(1)\tilde{u}-2\lambda \tilde{\chi}\tilde{u}(1).
\end{aligned}
\end{equation}
Now by recalling that $(\tilde{u}, \tilde{\chi})\in D(\mathcal{A}_e)$, using integration by parts and bounding one has
$$
\omega=-\gamma\tilde{u}(1)^2+\tilde{u}(1)\tilde{\chi}-\int_0^1 \tilde{u}_x^2 \leq -\gamma\tilde{u}(1)^2+\tilde{u}(1)\tilde{\chi}
$$
Plugging the above bound into \eqref{eq:Vdot}, integrating the term $\tilde{u}_{xx}$ and using again the fact that 
$(\tilde{u}, \tilde{\chi})\in D(\mathcal{A}_{e})$ one gets
\begin{equation}
\begin{aligned}
 \dot{V}(\tilde{z})\leq &-2k \int_0^1\tilde{u}^2-2\lambda\gamma\tilde{u}(1)^2+\varepsilon \int_0^1 k\tilde{u}\tilde{\chi}\\
 &+\varepsilon\rho \int_0^1 \tilde{u}(1)\tilde{u}-\lambda\varepsilon \tilde{\chi}^2+\lambda\varepsilon\gamma\tilde{u}(1)\tilde{\chi}.
\end{aligned}
\label{eq:Vdot_bound0}
\end{equation}
By defining 
$$
\Psi(\varepsilon)\coloneqq\begin{bmatrix}
    -2k &\frac{\varepsilon k}{2}&\frac{\varepsilon\rho}{2}\\
    \star&-\lambda\varepsilon&\frac{\lambda\varepsilon\gamma}{2}\\
    \star&\star&-2\gamma\lambda
\end{bmatrix}
$$
\eqref{eq:Vdot_bound0} can be equivalently rewritten as
\begin{equation}
\label{eq:Vdot2}
\dot{V}(\tilde{z})\leq\int_0^1 \begin{bmatrix}
    \tilde{u}\\
    \tilde{\chi}\\
    \tilde{u}(1)
\end{bmatrix}^\top \Psi(\varepsilon) \begin{bmatrix}
    \tilde{u}\\
    \tilde{\chi}\\
    \tilde{u}(1)
\end{bmatrix}.
\end{equation}
Next we show that there exists $\varepsilon^\star>0$ such that for all $\varepsilon\in (0, \varepsilon^\star]$,
$\Psi(\varepsilon)\prec 0$. To this end, notice that by Schur's complement lemma $\Psi(\varepsilon)\prec 0$ if and only if 
$$
\Upsilon(\varepsilon)\coloneqq \underbrace{\begin{bmatrix}
    -2k &\frac{\varepsilon k}{2}\\
    \star&-\lambda\varepsilon
\end{bmatrix}}_{\mathcal{L}_\Upsilon(\varepsilon)}+\frac{\varepsilon^2}{8\lambda\gamma}\begin{bmatrix}
   \rho^2&\lambda\gamma\rho\\
    \star&\lambda^2\gamma^2
\end{bmatrix}\prec 0.
$$
Therefore, by using Weyl's inequality (see, e.g., \cite[Theorem 8.4.11]{bernstein2009matrix}) 
$$
\begin{aligned}
\lambda_{\max}(\Upsilon(\varepsilon))&=\lambda_{\max}(\Upsilon(\varepsilon)+\mathcal{L}_\Upsilon(\varepsilon)-\mathcal{L}_\Upsilon(\varepsilon))\\
&\leq \lambda_{\max}((\Upsilon(\varepsilon))+ \lambda_{\max}(\Upsilon(\varepsilon)-\mathcal{L}_\Upsilon(\varepsilon)).
\end{aligned}
$$
The latter, by recalling that, due to symmetry, $\Vert \Upsilon(\varepsilon)-\mathcal{L}_\Upsilon(\varepsilon)\Vert_2=\vert \lambda_{\max}(\Upsilon(\varepsilon)-\mathcal{L}_\Upsilon(\varepsilon)\vert$ gives
\begin{equation}
\begin{aligned}
\lambda_{\max}(\Upsilon(\varepsilon))\leq \lambda_{\max}((\mathcal{L}_\Upsilon(\varepsilon))+ \Vert \Upsilon(\varepsilon)-\mathcal{L}_\Upsilon(\varepsilon)\Vert_2
\end{aligned}
\label{eq:Upsilon_neg2}
\end{equation}
Simple calculations show that
\begin{equation}
\label{eq:Upsilon_neg}
\mathcal{L}_\Upsilon(\varepsilon)\prec 0,\quad \forall \varepsilon\in \left(0, \frac{8\lambda}{k}\right).
\end{equation}
Now, let 
\begin{equation}
\label{eq:a_b}
0<a<b<\frac{8\lambda}{k}
\end{equation}
and define
\begin{equation}
M\coloneqq\max_{\varepsilon\in [a, b]}\lambda_{\max}(\mathcal{L}_\Upsilon(\varepsilon)).
\label{eq:M_scalar}
\end{equation}
Notice that $M$ is well-defined due to $\mathcal{L}_\Upsilon$ being continuous on $\varepsilon$ and that $M<0$ thanks to \eqref{eq:Upsilon_neg} and the selection of $a$ and $b$ in \eqref{eq:a_b}. Pick
\begin{equation}
\label{eq:delta}
\delta \in \left(0, \frac{1}{b}\vert M\vert\right).
\end{equation}
Select $\varepsilon^\prime>0$ such that 
\begin{equation}
\Vert \Upsilon(\varepsilon)-\mathcal{L}_\Upsilon(\varepsilon)\Vert_2\leq \delta \varepsilon\quad\forall\varepsilon\in (0, \varepsilon^\prime).
\label{eq:NormBound}
\end{equation}
This is always possible since $\varepsilon\mapsto\Upsilon(\varepsilon)$ is Fréchet differentiable and by construction
$$
\Upsilon(0)+d\Upsilon(0)\varepsilon=\mathcal{L}_\Upsilon(\varepsilon),\quad\forall \varepsilon\in\R.
$$
Then, by combining \eqref{eq:Upsilon_neg2} and \eqref{eq:NormBound} one has 
\begin{equation}
\lambda_{\max}(\Upsilon(\varepsilon))\leq \lambda_{\max}((\mathcal{L}_\Upsilon(\varepsilon))+\delta \varepsilon\quad\forall\varepsilon\in (0, \varepsilon^\prime)
\label{eq:Bound_var_ep_prime}
\end{equation}
Now pick $\varepsilon^\star\in (a, \min\{b,\varepsilon^\prime\})$. Then, since $\varepsilon^\star\in (0, \varepsilon)$, from \eqref{eq:Bound_var_ep_prime} one gets
$
\lambda_{\max}(\Upsilon(\varepsilon^\star))\leq \lambda_{\max}((\mathcal{L}_\Upsilon(\varepsilon^\star))+\delta b
$
Hence, from the definition of $M$ in \eqref{eq:M_scalar}, the latter yields 
$
\lambda_{\max}(\Upsilon(\varepsilon^\star))\leq M+\delta b
$
which by using the selection of $\delta$ in \eqref{eq:delta} gives
$
\lambda_{\max}(\Upsilon(\varepsilon^\star))< M+\vert M\vert=0.
$

To conclude the proof, let  $\varepsilon^\star>0$ such that $\Psi(\varepsilon^\star)\prec 0$. Pick 
$\varepsilon\in \left(0, \min\left\{\varepsilon^\star,\sqrt{\frac{\lambda}{\rho}}\right\} \right)$. Then, from \eqref{eq:Vdot2} one has
\begin{equation}
  \dot{V}(\tilde{z})\leq -\vert  \lambda_{\max}(\Psi(\varepsilon))\vert  ( \langle \tilde{u},\tilde{u}\rangle+\tilde{\chi}^2)
\end{equation}
which by using 
\eqref{eq:sandwich} and the definition of the inner product $\langle\cdot, \cdot \rangle_{\mathcal{X}}$ yields
\begin{equation}
\label{eq:Vdot_final}
  \dot{V}(\tilde{z})\leq -\vert  \lambda_{\max}(\Psi(\varepsilon))\vert\max\left\{1, \frac{\rho}{\lambda}\right\}V(\tilde{z}).
\end{equation}
This concludes the proof.
\end{proof}
\subsection{Recovery of the unknown parameter}
The steady-state $\bar\chi$ can be expressed as $\bar{\chi}=g_{\lambda}(k)y_r$, where $g_{\lambda}:[0,\infty)\rightarrow[0,\infty)$ is the function defined by
\begin{equation}\label{eq:function_g}
g_{\lambda}(k)=\sqrt{\frac{k}\lambda}\tanh\left(\sqrt{\frac{k}\lambda}\right)
\end{equation}
The following technical lemma establishes some properties of the function $g_{\lambda}(k)$ that will be helpful to devise an algorithm for the recovery of the unknown parameter $k$.
\begin{lemma}\label{lem:function_g}
{\it For any fixed $\lambda>0$, the function $g_{\lambda}(\cdot)$ is a class $\mathcal{K}^{\infty}$ function, so that the inverse function $g_{\lambda}^{-1}(\cdot)$ exists and is well defined over $[0,\infty)$. Moreover, $g_{\lambda}(\cdot)$ is continuosuly differentiable in $(0,\infty)$ with $g_{\lambda}'(\cdot)>0$.}
\end{lemma}\smallskip
Based on the previous result, we can give an algorithm for the recovery of $k$ from the estimate $\chi$, which  follows as a direct corollary from Proposition~\ref{prop:GES}.
\begin{corollary}
{\it Consider the error system \eqref{eq:error_dyn}, with $y_r\neq 0.$ Define the estimated parameter
\begin{equation}\label{eq:khat}
\hat{k}(t):=g_{\lambda}^{-1}\left(\max\left\{0,\frac{{\chi}(t)}{y_r}\right\}\right)
\end{equation}
Then one has
$$
\lim_{t\rightarrow+\infty}|k-\hat{k}(t)|=0
$$}
\end{corollary}
\begin{proof}
Thanks to the convergence of $\chi(t)$ to $\bar{\chi}$, by construction there exists $\bar{t}\geq0$ such that $\chi(t)/\bar{\chi}>0$ for any $t\geq\bar{t}$, which implies also that
$\hat{k}(t)=g_{\lambda}^{-1}(\chi(t)/y_r)$ for $t\geq \bar{t}$. Recalling that the function $g_{\lambda}^{-1}(\cdot)$ is continuous, one has
$$
\lim_{t\rightarrow+\infty}\hat{k}(t)=g_{\lambda}^{-1}\left(\frac{\bar{\chi}}{y_r}\right)=g_{\lambda}^{-1}(g_{\lambda}(k))=k
$$
and this proves the claim.
\end{proof}
It is worth stressing that, even though the inverse function $g^{-1}_\lambda(k)$ is well-defined, no closed-form expression can be found and therefore \eqref{eq:khat} needs to be implemented via numerical evaluation.
\section{Joint state and parameter estimation}\label{sec:joint}
The estimation procedure exemplified earlier can be also used in combination with an infinite-dimensional observer to retrieve the full state of the system $u(x,t)$ from the boundary measurement $y(t)$. To this end, let us design an observer by copying the original dynamics of the system and implementing an output injection at the right boundary, with a gain $\alpha>0$. Clearly, as the true value $k$ is unknown, the observer dynamics can only be defined through the estimated value $\hat{k}(t)$ defined by \eqref{eq:khat}. We have then
\begin{equation}\label{eq:observer}
\begin{array}{rcl}
\hat{u}_t(x,t)&=&\lambda \hat{u}_{xx}(x,t)-\hat{k}(t) \hat{u}(x,t),\,\, x\in[0,1],\, t\geq 0\\
\hat{u}_x(0,t)&=&0\\
\hat{u}_x(1,t)&=&v(t)+\alpha(y(t)-\hat{u}(1,t))
\end{array}
\end{equation}
Let us define the estimation error $\eta(x,t)\coloneqq u(x,t)-\hat{u}(x,t)$, whose dynamics is governed by the PDE
\begin{equation}\label{eq:err_obs}
\begin{array}{rcl}
\eta_t(x,t)&=&\lambda 
\eta_{xx}(x,t)-ku(x,t)+\hat{k}(t)\hat{u}(x,t)\\
&=&\lambda 
\eta_{xx}(x,t)-\hat{k}(t)\eta(x,t)-(k-\hat{k}(t))u(x,t)\smallskip\\
\eta_x(0,t)&=&0\\
\eta_x(1,t)&=&-\alpha\eta(1,t)
\end{array}
\end{equation}
Let us rewrite \eqref{eq:err_obs} as the following abstract semilinear dynamical system 
\begin{equation}\label{eq:abstracterr}
\begin{aligned}
&\left[\begin{array}{l}
\dot{\eta}\\
\dot{\tilde{u}}\\
\dot{\tilde{\chi}}
\end{array}\right]=\widehat{A}
\begin{bmatrix}
\eta\\
\tilde{u}\\
\tilde{\chi}
\end{bmatrix}+f(\tilde{\chi}, \tilde{u},\eta)
  \end{aligned}
\end{equation}
where
$$
\begin{aligned}
&\widehat{A}\coloneqq\begin{bmatrix}
A_\eta&0\\
0&A_e\end{bmatrix}\\
&D(A_\eta)\coloneqq \left\{\eta\in\mathcal{H}^2(0,1)\colon \eta_x(0)=0, \eta_x(1)=-\alpha\eta(1)\right\}\\
&A_\eta \eta\coloneqq \lambda \eta_{xx}\\
&f(\tilde{\chi}, \tilde{u},\eta)\coloneqq\begin{bmatrix}
-\Psi(\tilde{\chi})\eta-(k-\Psi(\tilde{\chi}))(\tilde{u}+u_r)\\
0
\end{bmatrix}
\end{aligned}
$$
and $A_e$ is defined in \eqref{eq:Ae_op}. The following result establishes well-posedness of system \eqref{eq:abstracterr}.
\begin{proposition}
System \eqref{eq:abstracterr} is well-posed in the sense of Definition~\ref{def:well-posed}. In particular, $f$ is locally Lipschitz continuous and the operator $\widehat{A}$ is the infinitesimal generator of a strongly continuous semigroup on the Hilbert space $\mathcal{L}_2(0,1)\oplus\mathcal{L}_2(0,1)\oplus\R$.
\end{proposition}
The proof of the result is wholly similar to the proof of Proposition~\ref{prop:wellposed} (the local Lipschitzness of $f$ follows directly from the definition of $f$ itself). Therefore, due to lack of space, it is omitted. 

Let us now consider the total error system, denoted by $\Sigma_{\mathrm{tot}}=(\tilde{u},\tilde{\chi},\eta)$, and defined as the joint system given by \eqref{eq:error_dyn} together with \eqref{eq:abstracterr}. 
Using the identity $u=\tilde{u}+u_r$ and defining\footnote{Observe that, by construction, we have $\Psi(0)=k$.}
$$
\hat{k}=\Psi(\tilde{\chi})\coloneqq g_{\lambda}^{-1}\left(\max\left\{0,\frac{{\tilde\chi}}{y_r}+\frac{{\bar\chi}}{y_r}\right\}\right)
$$
we can observe that the total error system $\Sigma_{\mathrm{tot}}$ is structured as a cascade nonlinear feedback of $(\tilde{u},\tilde{\chi})$ with $\eta$. In particular, the dynamics of $\eta$ in \eqref{eq:abstracterr} can be alternatively rewritten as
\begin{equation}\label{eq:psi-eta}
\dot{\eta}=A_\eta \eta-\Psi(\tilde{\chi})\eta-(k-\Psi(\tilde{\chi}))(\tilde{u}+u_r)
\end{equation}
It is worth stressing that $u_r\in\mathcal{L}^2(0,1)$ is considered to be arbitrary but fixed here, resulting in the presence of  affine terms in the dynamics \eqref{eq:psi-eta}.
Before proceeding with the analysis of the stability of the total error system, let us exploit the next fact, following from~Proposition~\ref{prop:GES} as a corollary.
\begin{fact}\label{fact:bounds}
{\it In light of the global exponential stability of $\tilde{u}(x,t)$ and $\tilde{\chi}(t)$, for any pair of initial conditions $\xi_0:=(\tilde{u}(x,0),\tilde{\chi}(0))\in\mathcal{L}^2(0,1)\times \mathbb{R}$ and any $\nu>0$, there exists $T_0=T_0(\|\xi_0\|,\nu)>0$ such that for any $t\geq T_0$ one has}
$$ 
\begin{array}{rl}
\|\tilde{u}(x,t)+u_r(x)\|_{\mathcal{L}^2}&\leq \nu+ \|u_r(x)\|_{\mathcal{L}^2}\smallskip\\
|k-\Psi(\tilde{\chi}(t))|&\leq \nu 
\end{array}
$$
\end{fact}
We are now ready to state the following convergence result.
\begin{theorem}\label{th:joint}
{\it The origin of the total error system $\Sigma_{\mathrm{tot}}$, governed by \eqref{eq:error_dyns} and \eqref{eq:psi-eta}, is GAS and LES.}
\end{theorem}
\begin{proof}
The proof is divided in two steps. First we prove that the origin of the total system is globally attractive, and that the system trajectories reach in finite-time a bounded and closed set containing the origin and having arbitrarily small size. Then the local exponential stability of the origin is established which, building on the global attractivity condition, entails global asymptotic stability as well.\\
To study the convergence of the total error system, let us begin by introducing the Lyapunov functional
$$
W(\eta):=\frac12\int_0^1\eta(x)^2dx
$$
Evaluating the derivative along the solutions of \eqref{eq:psi-eta}, upon integration by parts and the application of Poincaré-Wirtinger inequality, yields
\begin{equation}\label{eq:Wdot}
\begin{array}{rcl}
\dot{W}(\eta)&=&-(\tilde{c}_0+2\Psi(\tilde{\chi}))W(\eta)\\
&&\displaystyle-(k-\Psi(\tilde{\chi}))\int_0^1\eta(x)(\tilde{u}(x)+u_r(x))dx\\
&\leq&\!\!\!\!\!\displaystyle-\tilde{c}_0W(\eta)+|k\!-\!\Psi(\tilde{\chi})|\!\int_{0}^1\!\!|\eta(x)(\tilde{u}(x)\!+\!u_r(x))|dx
\end{array}
\end{equation}
where $\tilde{c}_0>0$ depends on $\lambda, \alpha$ and the Poincaré-Wirtinger constant $c_{\mathrm{P}}$. Now, invoking H\"older's inequality, the previous condition can be further manipulated as
\begin{equation}\label{eq:ISS}
\dot{W}(\eta)\leq-\tilde{c}_0 W(\eta)+\sqrt{2}W(\eta)^{\frac12}|k-\Psi(\tilde{\chi})|\|\tilde{u}+u_r\|_{\mathcal{L}^2}
\end{equation}
Let us now show that, for any $\varsigma>0$, there exists a finite number $T_{\varsigma}>0$ such that solutions are confined in the set
$$
\mathcal{B}_{\varsigma}\coloneqq\left\{\eta\in\mathcal{L}^2(0,1): W(\eta)=\frac12\|\eta\|_{\mathcal{L}^2}^2\leq \varsigma\right\}
$$
for $t\geq T_{\varsigma}$, 
 thereby proving that $\{\eta(x)=0\}$ is globally attractive. To this end, fix $\varsigma>0$. By exploiting Fact~\ref{fact:bounds} and observing that, from inequality \eqref{eq:ISS}, $W(\eta)$ is qualified as ISS-Lyapunov function for \eqref{eq:psi-eta}, we have
$$
\begin{array}{rcl}
\dot{W}(\eta)&\leq&-\tilde{c}_0W(\eta)+\sqrt{2}W(\eta)^{\frac12}(\nu^2+\nu\|u_r(x)\|_{\mathcal{L}^2})\smallskip\\
&\leq&-\displaystyle\frac{\tilde{c}_0}2W(\eta)
\end{array}
$$
for any $t\geq T_0(\xi_0,\nu)$ and as long as $\eta(t)$ satisfies
$$
W(\eta(t))\geq \frac{8}{\tilde{c}_0^2}(\nu^2+\nu\|u_r\|_{\mathcal{L}^2})^2.
$$
In particular, the latter inequality guarantees that $\eta$ is bounded.
Moreover, since $\nu$ can be chosen arbitrarily, by picking $\nu<\nu_{\varsigma}$ with
$$
\frac8{\tilde{c}_0^2}(\nu_\varsigma^2+\nu_\varsigma\|u_r(x)\|_{\mathcal{L}^2})^2=\varsigma,
$$
the claimed attractivity property holds uniformly with $T_{\varsigma}=T_0(\|\xi_0\|,\nu_\varsigma)$. Combining such property with the stability of the dynamics~\eqref{eq:error_dyns}, we can easily infer that, for any closed and bounded set $\mathcal{E}\subset\mathcal{L}^2\times\mathbb{R}$ and for any given $\varsigma>0$, the convergence of $(\tilde{z},\eta)$ onto the set  $\mathcal{E}\times\mathcal{B}_{\varsigma}$ occurs in finite-time. This shows indeed that the  equilibrium $\{\tilde{z},\eta=(0,0)\}$ of the total system is globally attractive.\\
To prove LES, consider now the composite Lyapunov functional 
$$
U(\tilde{z},\eta):=V(\tilde{z})+\vartheta W(\eta)
$$
with $V(\tilde{z})$ defined as in \eqref{eq:Lyapunov} and $\vartheta>0$ to be selected. From the Lyapunov analysis carried out in the proof of Proposition~\ref{prop:GES}, we can infer that
$$
\dot{V}(\tilde{z})\leq - \phi_1  V(\tilde{z})-\phi_2|\tilde{\chi}|^2
$$
for some positive constants\footnote{For example one can pick $\phi_1=\sigma\vert  \lambda_{\max}(\Psi(\varepsilon))\vert\max\left\{1, \frac{\rho}{\lambda}\right\}$ and $\phi_2=(1-\sigma)\vert  \lambda_{\max}(\Psi(\varepsilon))\vert\max\left\{1, \frac{\rho}{\lambda}\right\}$ for some $\sigma\in(0,1)$.} $\phi_1,\phi_2>0$. 
For what concerns the handling of $\dot{W}(\eta)$, consider again \eqref{eq:Wdot} and apply Young's inequality to get 
$$
\dot{W}(\eta)\leq -\left(\tilde{c}_0-\delta\right)W(\eta)+\frac{1}{2\delta}|k-\Psi(\tilde{\chi})|^2\|\tilde{u}(x)+u_r(x)\|_{\mathcal{L}^2}^2
$$
with $0<\delta<\tilde{c}_0$. Using again Fact~\ref{fact:bounds}, there exists an open region of initial conditions ${\mathcal{E}}$  such that, if $\tilde{z}(0)\in\mathcal{E}$, then
$$
\|u(x)\|_{\mathcal{L}^2}^2\leq (1+\|u_r(x)\|_{\mathcal{L}^2})^2=:c_3\ \forall t\geq0.
$$
Furthermore, using a Taylor expansion around $\tilde{\chi}=0$ and applying the formula for the inverse function derivative, the term $k-\Psi(\tilde{\chi})$ can be expressed as
$$
\begin{array}{rl}
k-\Psi(\tilde{\chi})&\displaystyle=k-g_{\lambda}^{-1}\left(\frac{\tilde{\chi}}{y_r}+g_{\lambda}(k)\right)\smallskip\\
&=\displaystyle-\frac{1}{g_{\lambda}'(k)}\frac{\tilde{\chi}}{y_r}+o(\tilde{\chi})
\end{array}
$$
where $g_{\lambda}(k)$ is defined as in \eqref{eq:function_g} and $o(\tilde{\chi})$ indicates higher-order terms. The above expression is always well-defined since, by definition, one has $g'_\lambda(k)>0$ for any $k>0$. Now, an open subset of initial conditions $\mathcal{E}_1\subseteq\mathcal{E}$ can always be chosen such that if $\tilde{z}\in\mathcal{E}_1$, then
$$
|o(\tilde{\chi})|\leq \left|\frac{\tilde{\chi}}{g_{\lambda}'(k)y_r}\right|\ \forall t\geq 0 
$$thus providing the estimate
$
|k-\Psi(\tilde{\chi})|^2\leq c_4 |\tilde\chi|^2\quad \forall t\geq 0
$
with $c_4:=4/(|g_{\lambda}'(k)y_r|)^2$. Now, putting all pieces together, we have shown that, for $(\tilde{z},\eta)\in\mathcal{E}_1\times \mathcal{L}^2(0,1)$, the following condition is met
$$
\begin{array}{rl}
\dot{U}(\tilde{z},\eta)&=\dot{V}(\tilde{z})+\vartheta \dot{W}(\eta) \smallskip\\
&\leq\displaystyle{- \phi_1  V(\tilde{z})}-\vartheta (\tilde{c}_0-\delta)W(\eta)\smallskip\\
&\displaystyle-\left({\phi_2}-\vartheta\frac{c_3c_4}{2\delta}\right)|\tilde{\chi}|^2
\end{array}
$$
Finally, by picking $\vartheta<(2\delta\phi_2)/(c_3c_4),
$ it is straighforward to verify that
$$
\dot{U}(\tilde{z},\eta)\leq -c_5 U(\tilde{z},\eta)
$$
with $c_5:=\min\{\phi_1,(\tilde{c}_0-\delta)\}$,
and this inequality implies the claimed local exponential stability condition. By the global attractivity property that has been established earlier, for any initial condition $\tilde{z}(0)$, there exists a finite time $T_1>0$ such that $\tilde{z}(t)\in\mathcal{E}_1$ for any $t\geq 0$, so that also GAS of the equilibrium $\{(\tilde{z},\eta)=(0,0)\}$ has been proved. This concludes the proof.
\end{proof}

\section{Numerical example and simulations}\label{sec:ex}
\begin{figure}[h!]
\hspace{-0.8cm}
\centering
\includegraphics[scale = 0.52]{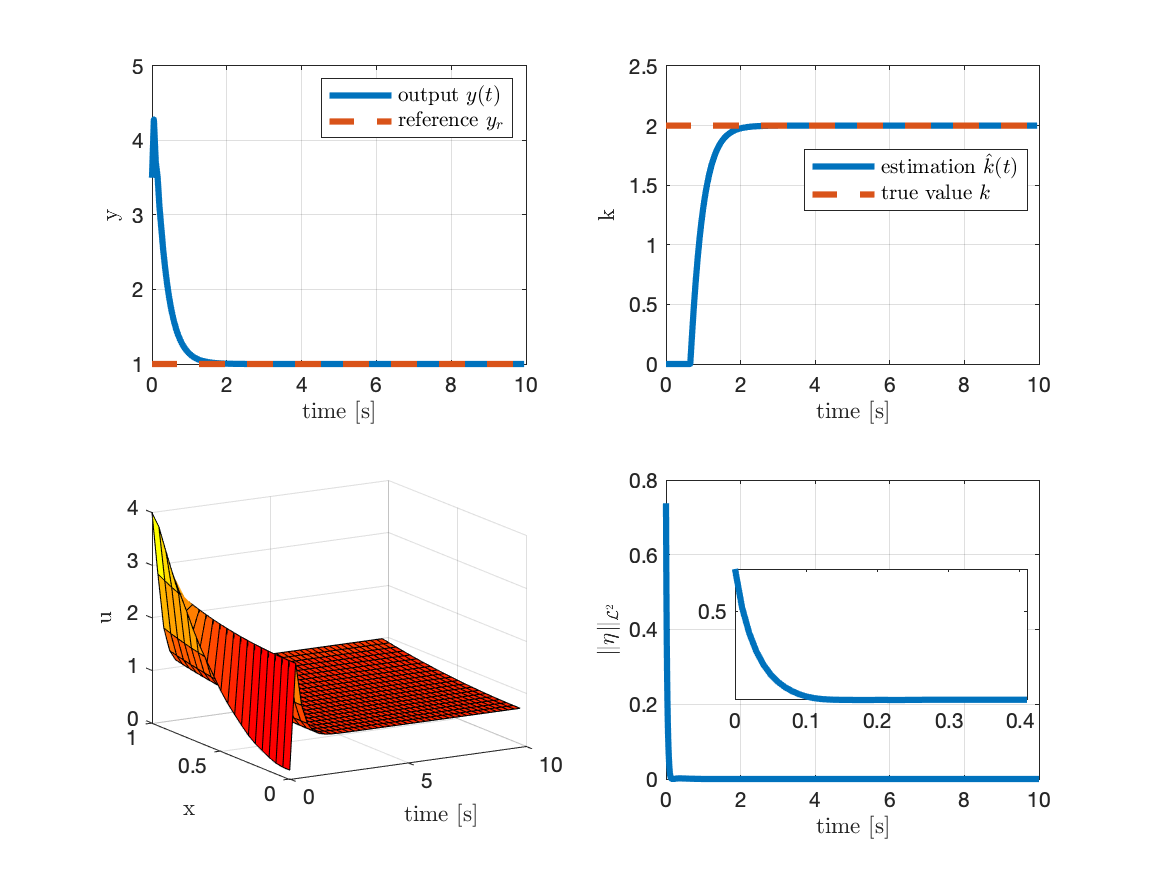}
\caption{Estimation results}\label{fig:example}
\end{figure}

Consider the reaction-diffusion equation \eqref{eq:sys}, with coefficients
$
\lambda=3,\ k=2
$
and initial condition
$
u_0(x)=7x^2/4
$.
Selecting the output reference $y_r=1$, we have implemented the estimation procedure detailed in Section~\ref{sec:main}. The adaptation gains in \eqref{eq:error_dyns} have been chosen as 
$
\gamma=2,\ \rho=4.5.
$
The simulations have been performed in \texttt{Matlab}, using a finite-dimensional approximation of the partial differential equations based on spectral decomposition \cite{kato2013perturbation}.
The results are shown in Figure~\ref{fig:example}. In particular, in the top-left plot, we can observe the output tracking performances of the control, whereas, in the top-right plot, we can appreciate the convergence of the estimation $\hat{k}(t)$ to the actual value of the unknown reaction coefficient $k$. The bottom-left plot shows the graph of the system solution. To further illustrate the capabilities of the proposed approach, the joint state and parameter estimation procedure detailed in Section~\ref{sec:joint} has been tested. In particular an observer with the structure \eqref{eq:observer} has been implemented. The bottom-right plot of Figure~\ref{fig:example} displays the obtained results, which are characterized by a very fast convergence of the estimation error. {In particular we can observe that the transient of the estimation error is significantly shorter than the one of the recovery of the reaction coefficient. This feature suggests that the proposed set-point regulation control is inherently robust to the uncertainty in the reaction coefficient.}
\section{Conclusions and extensions}\label{sec:conc}
This paper deals with a reaction-diffusion equation with collocated boundary control and unknown reaction coefficient. Avoiding the use of classical adaptive control tools, we proposed here a novel estimation technique for the unknown coefficient based on the use of set-point regulation controller. In particular, hinging on the particular structure of the steady-state solution, the value of the reaction coefficient can be recovered by nonlinear mapping inversion. The obtained results are then shown to be pivotal for the design of a state observer for the uncertain PDE based on boundary output only, which allow for the simultaneous estimation of the reaction coefficient and the full state of the parabolic equation, owing to the cascade structure of the total error system.

We are currently working on the extension of the proposed adaptive estimation technique to reaction-advection-diffusion equations and, more in general,  to systems of linear parabolic equations with uncertain coefficients. The use of time-varying output references $y_r(t)$ could also be object of future investigations.
\bibliographystyle{ieeetran}
\bibliography{biblio}
\end{document}